\documentclass [11pt, article, oneside, a4paper]{memoir}

\usepackage{amsmath,amsfonts,amssymb,amsthm,paralist,graphicx,xargs}
\usepackage[all]{xy}
\usepackage[backend=biber,maxbibnames=15]{biblatex}
\usepackage[T1]{fontenc}
\usepackage{amsmath,amsfonts,lmodern}
\usepackage[colorlinks=true,linkcolor=black,citecolor=blue]{hyperref}
\usepackage{paralist}
\usepackage{color}

\newcommand{\evM}[2]{\bigl[#1\bigr]_{#2}}
\newcommand{\evm}[2]{[ #1]_{#2}}
\newcommand{\up}[1]{\,\uparrow #1}
\newcommand{\Up}[1]{\,\Uparrow #1}
\newcommand{\bbR}{\mathbb{R}}

\newcommand{\To}{\Rightarrow}
\newcommand{\Tos}{[[\,\Rightarrow\,\!]]}
\newcommand{\tq}{\;:\;}

\DeclareMathOperator{\Id}{Id}

\newcommand{\D}{\mathcal{D}}
\renewcommand{\H}{\mathsf{H}}
\newcommand{\A}{\mathsf{A}}
\newcommand{\B}{\mathsf{B}}
\newcommand{\Y}{\mathsf{Y}}

\newcommand{\T}{\mathsf{T}}

\newcommand{\U}{\mathsf U}

\newcommand{\rtilde}{{\widetilde r}}
\newcommand{\elltilde}{{\widetilde \ell}}

\newcommand{\Obar}{\overline O}
\newcommand{\Tbar}{\overline{T}}

\newcommand{\Ttilde}{\widetilde{T}}
\newcommand{\Ytilde}{\widetilde{\Y}}

\newcommand{\initial}[1]{\iota(#1)}
\newcommand{\final}[1]{\kappa(#1)}
\newcommand{\length}[1]{\ell(#1)}
\newcommand{\w}{\mathsf{w}}
\newcommand{\X}{\mathsf{X}}
\newcommand{\ie}{\textsl{i.e.}}
\newcommand{\eof}{\ \rule{.2em}{1em}}

\newcommand{\addpoint}[1]{#1\ ---\ }

\setsecnumdepth{subsubsection}

\setsecnumformat{\csname the#1\endcsname---}
\setsubsubsechook{\setsecnumformat{{\normalfont\sffamily\csname
the##1\endcsname\ }}}
\setsubsechook{\setsecnumformat{\csname the##1\endcsname\ ---\ }}
\setsechook{\setsecnumformat{\csname the##1\endcsname---}}

\setsecheadstyle{\centering\Large\bfseries\sffamily}
\setsubsecheadstyle{\large\bfseries\sffamily}
\setsubsubsecheadstyle{\bfseries\itshape\addpoint}
\setsubsubsecindent{1em}
\setbeforesubsubsecskip{.5em plus .2em minus -.1em}
\setaftersubsubsecskip{-0em}
\setsubparaheadstyle{\itshape}
\newtheoremstyle{thm}%
{1.5ex plus .3ex minus .1ex}%
{1ex plus .3ex minus .1ex}%
{\itshape}%
{}%
{\sffamily}%
{---}%
{0em}%
{$\bullet$\hbox{\ }#1\hbox{\ }#2}%

\theoremstyle{thm}

\newtheorem{definition}{Definition}[section]

\newtheorem{theorem}[definition]{Theorem}

\newtheorem{lemma}[definition]{Lemma}
\newtheorem{proposition}[definition]{Proposition}
\newtheorem{corollary}[definition]{Corollary}

\newtheoremstyle{note}%
{1ex plus .3ex minus .1ex}%
{1ex plus .3ex minus .1ex}%
{}%
{}%
{\itshape}%
{.}%
{1em}%
{}%
\theoremstyle{note}
\newtheorem{example}[definition]{Example}
\newtheorem{remark}[definition]{Remark}

\begin{document}
\begin{center}
  \Large\textbf{Convergence of distributions on paths\footnote{Published in H. Fernau and K. Jansen, editors, \emph{Fundamentals of Computation Theory} (Proceedings of FCT 2023, Tier, Germany), LNCS 14292, Springer, pp.~1--15, 2023.}}\par
  \bigskip
  \large
  Samy Abbes\\
  \small Université Paris Cité, CNRS\\ IRIF, F-75013 Paris, France\\
\texttt{abbes@irif.fr}
\end{center}

\begin{abstract}
We study the convergence of distributions on finite paths of weighted digraphs, namely the family of Boltzmann distributions and the sequence of uniform distributions. Targeting applications to the convergence of distributions on paths, we revisit some known results from \emph{reducible} nonnegative matrix theory and obtain new ones, with a systematic use of tools from analytic combinatorics. In several fields of mathematics, computer science and system theory, including concurreny theory, one frequently faces non strongly connected weighted digraphs encoding the elements of combinatorial structures of interest; this motivates our study.

\medskip\small\textbf{Keywords:} Weighted digraph; Reducible nonnegative matrix; Uniform measure; Boltzmann measure
\end{abstract}

\section{Introduction}
\label{sec:introduction}

\paragraph{Motivations.}

Given a weighted digraph, it is standard to consider for each vertex $x$ and for each integer $k\geq0$ the finite probability distribution $\mu_k$ which gives to each path of length $k$ and starting from $x$ the probability proportional to its multiplicative weight. If the underlying digraph is strongly connected and aperiodic, the classical results of Perron-Frobenius theory for primitive matrices show that the sequence $(\mu_k)_{k\geq0}$ converges weakly toward a probability measure on the space of infinite paths starting from~$x$. The sequence of random vertices visited by such a uniform infinite path is a Markov chain, which transition kernel is derived from the weights on the edges and from the Perron eigenvector of the adjacency matrix of the weighted digraph~\cite{parry64}.

The same question of convergence is less standard in the case of a general, non strongly connected weighted digraph. It is however the actual framework for a variety of situations in mathematics, in computer science and in system theory. Indeed, the elements of many combinatorial or algebraic structures of interest can be represented by the paths of some finite digraph encoding the constraints of the ``normal form'' of elements, in a broad sense. Even when the initial structure is irreducible in some sense, it might very well be the case that the digraph itself being not strongly connected. A first example is the digraph of simple elements of a braid group or of a braid monoid, or more generally of a Garside group or monoid~\cite{charney95,dehornoy14,abbes19}; this digraph is never strongly connected. But also in automatic group theory: quoting~\cite{calegari10}, ``a graph parameterizing a combing of a hyperbolic group may typically fail to be recurrent''. In system theory, the digraph of states-and-cliques introduced by the author to encode the trajectories of a concurrent system~\cite{abbes22} is also not strongly connected in general, even for an irreducible concurrent system. 

Understanding the asymptotic behavior of ``typical elements of large length'' is of great interest at some point, either for random generation purposes, or aiming for the probabilistic verification or the performance evaluation of systems, to name a few applications. This motivates the study of the weak limit of the uniform distributions on paths for a general weighted digraph, which is the main topic of this paper.

\paragraph{Framework and contributions.}
\label{sec:fram-contr}

Let $W=(V,\w)$ be a weighted digraph, where $\w:V\times V\to\bbR_{\geq0}$ is a non negative weight function, positive on edges. Let $F$ be the adjacency matrix of~$W$, that we  assume to be of positive spectral radius~$\rho$.

Let $\H(z)$ be the generating function with matrix coefficients defined by $\H(z)=\sum_{k\geq0}F^kz^k$. The power series $\H(z)$, that we call the \emph{growth matrix} of~$W$, can also be seen as a matrix of generating functions with a well known combinatorial interpretation, which is recalled later in the paper.
%
%
It is standard knowledge that all the generating functions $\evm{\H(z)}{x,y}$ are rational series. From this we can deduce the existence and uniqueness of a triple $(\lambda,h,\Theta)$ where $\lambda$ is a positive real, $h$~is a positive integer and $\Theta$ is a nonnegative and non zero square matrix indexed by $V\times V$, that we call the \emph{residual matrix of\/~$W$}---it can indeed be interpreted as the residue at its positive smallest singularity of the function~$\H(z)$---, and such that:
\begin{gather}
  \label{eq:2}
  \lim_{s\to\lambda}(1-\lambda^{-1}s)^h\,\H(s)=\Theta\qquad s\in(0,\lambda)
\end{gather}

For instance if $W$ is strongly connected and aperiodic, and so $F$ is primitive, then it is well known that there is projector matrix $\Pi$ of rank $1$ and a square matrix $R$ of spectral radius lower than $1$ such that:
\begin{align}
  \label{eq:3}
F&=\rho(\Pi+R)&  \Pi\cdot R&=0&R\cdot\Pi&=0  
\end{align}
In this case: $\lambda=\rho^{-1}$, $\Theta=\Pi$, $h=1$, and the elements $\evm\Theta{x,y}$ are all positive.

For the less standard situation where $W$ is not strongly connected, a contribution of this paper is to provide a recursive way of computing the residual matrix~$\Theta$ and to characterize the pairs $(x,y)$ such that $\evm\Theta{x,y}>0$ (Th.~\ref{thr:1} and Sec.~\ref{sec:comp-asympt-matr}). For this we use a simple formula (Lemma~\ref{lem:2}) to compute the rectangular blocks of $\H(z)$ outside its diagonal blocks---the latter are occupied by the corresponding growth matrices of the access classes of the digraph; despite its simplicity, this formula does not seem to have explicitly appeared in the literature.

We also show that the integer $h$ occurring in~(\ref{eq:2}) is the \emph{height} of~$W$, that is to say, the maximal length of chains of access equivalence classes of maximal spectral radius; in nonnegative matrix theory, access classes of maximal spectral radius are called \emph{basic}. This is consistent with the known result saying that the dimension of the generalized eigenspace associated to $\rho$ is the height of the digraph~\cite{rothblum14}, a result that we also recover. This yields precise information on the growth of coefficients of the generating functions $\evm{\H(z)}{x,y}$.

We use these results to study in Sec.~\ref{sec:conv-distr} the weak convergence of certain probability distributions on paths, in particular the sequence of \emph{uniform distributions} and the family of \emph{Boltzmann distributions} (see Def.~\ref{def:3}). In the case of a general weighted digraph, we prove the convergence of the Boltzmann distributions toward a \emph{complete cocycle measure} on a sub-digraph (see Def.~\ref{def:4} and Th.~\ref{thr:7}). Under an aperiodicity assumption, the convergence also holds for the uniform distributions (Th.~\ref{thr:8}). 

The sub-digraph, support of the limit measure, has the property that its basic access classes coincide with its final access classes. We call such digraphs \emph{umbrella digraphs}, and we devote Sec.~\ref{sec:umbrella-digraphs} to their particular study. We give a decomposition of their adjacency matrix which generalizes the decomposition~(\ref{eq:3}) for primitive matrices (Th.~\ref{thr:2} and~\ref{thr:3}), and reobtain in this way some results from nonnegative matrix theory, for instance that umbrella digraphs are exactly those for which there exists a positive Perron eigenvector (Corollary~\ref{cor:1}). Our point of view is influenced by our motivation toward probability measures on the space of infinite paths, and we characterize umbrella digraphs by the existence of a complete cocycle measure (see Def.~\ref{def:4} and Th.~\ref{thr:6}).

\section{Preliminaries}
\label{sec:preliminaries}

\paragraph{\textbullet\bfseries Weighted digraphs, adjacency matrix, paths.}
\label{sec:preliminaries-1}

A \emph{weighted digraph}, or \emph{digraph} for short, is a pair $W=(V,\w)$ where $V$ is a finite set of \emph{vertices} and $\w:V\times V\to\bbR_{\geq0}$ is a nonnegative real valued function. The set of pairs $E=\{(x,y)\in V\times V\tq \w(x,y)>0\}$ is the set of \emph{edges} of~$W$. Given a bijection $\sigma:V\to\langle\nu\rangle$, where $\langle\nu\rangle=\{1,\ldots,\nu\}$, the function $\w$ identifies with the $\nu\times\nu$ nonnegative matrix $F$ defined by $\evm F{x,y}=\w(x,y)$. Changing the bijection $\sigma$ results in a simultaneous permutation of the lines and of the columns of~$F$.

A  \emph{finite path} of $W$ is a sequence $u=(x_i)_{0\leq i\leq k}$ of vertices such that  $(x_i,x_{i+1})\in E$ for all $i<k$. The \emph{initial} and \emph{final} vertices of $u$ are denoted $\initial u=x_0$ and $\final u=x_k$, and its \emph{length} is $\length u=k$. Denoting  by $O$ the set of finite paths, we introduce the following notations, for $x,y\in V$ and for $k$ any nonnegative integer:
\begin{align*}
  O_x&=\{u\in O\tq \initial u=x\}&O_{x,y}&=\{u\in O\tq \initial u=x\land\final u=y\}\\
  O_x(k)&=\{u\in O_x\tq \length u=k\}&O_{x,y}(k)&=\{u\in O_{x,y}\tq \length u=k\}
\end{align*}

 The \emph{infinite paths} are the sequences $\omega=(x_i)_{i\geq0}$ such that $(x_i,x_{i+1})\in E$ for all~$i\geq0$. The set $\Obar$ of paths, either finite or infinite, is equipped with the prefix ordering which we denote by~$\leq$. For every $u\in O$, $x\in V$  and $k\geq0$, we put:
 \begin{align*}
   \Omega_x&=\{\xi\in\Obar\tq\initial u=x\land\length u=\infty\}&   \Up u&=\{\xi\in\Obar\tq u\leq\xi\}\\
\up u&=\{\omega\in\Obar\setminus O\tq u\leq\omega\}&
\up^ku&=\{v\in O(k)\tq u\leq v\}
\end{align*}

\paragraph{\textbullet\bfseries Access relation, access equivalence classes.}
We write $x\To y$ to denote that the vertex $x$ has \emph{access} to the vertex~$y$, meaning that $O_{x,y}\neq\emptyset$. The relation $\To$ is then the transitive and reflexive closure of~$E$, seen as a binary relation on~$V$; it is thus a preordering relation on~$V$. We write $x\sim y$ to denote that $x$ and $y$ \emph{communicate}, meaning that $x\To y$ and $y\To x$, collapse equivalence relation of~$\To$. The equivalence classes of $\sim$ are the \emph{access classes} of~$W$. The set $\D$ of access classes is equipped with the partial ordering relation, still denoted~$\To$, induced by the preordering on vertices.

For each vertex~$x$, we define the sub-digraph:
\begin{align}
V(x)&=\{y\in V\tq x\To y\}
\end{align}

By convention, the bijection $\sigma:V\to\langle\nu\rangle$ will always be chosen in such a way that:\quad$\forall (x,y)\in V\times V\quad (x\To y)\implies \bigl(\sigma(x)\leq\sigma(y)\bigr)$. As a consequence, the adjacency matrix $F$ has the following block-triangular shape:
\begin{gather}
  \label{eq:5}
  F=\begin{pmatrix}
    F_1&\X&\dots&\X\\
    \vdots&\ddots&&\X\\
    0&\dots&0&F_p
  \end{pmatrix}
\end{gather}
where the diagonal blocks $F_1,\ldots,F_p$ are the adjacency matrices of the access classes, and the $\X$s represent rectangular nonnegative matrices.

\paragraph{\textbullet\bfseries Spectral radius, basic and final classes.}

By definition, the \emph{spectral radius} $\rho(W)$ of $W$ is the spectral radius $\rho(F)$ of its adjacency matrix~$F$, \ie, the largest modulus of its complex eigenvalues. It is apparent on~\eqref{eq:5} that $\rho(W)=\max_{D\in\D}\rho(D)$. In~\eqref{eq:5}, each matrix $F_i$ is either the $1\times 1$ block $[0]$ or an irreducible matrix. Hence it follows from Perron-Frobenius theory that $\rho(W)$ is itself an eigenvalue of~$F$.

By definition, an access class $D$ is \emph{basic} if $\rho(D)=\rho(W)$; and \emph{final} if $D$ is maximal in $(\D,\To)$.

For each vertex~$x$, we set:\quad$\gamma(x)=\rho\bigl(V(x)\bigr)$.

\paragraph{\textbullet\bfseries Analytic combinatorics, growth matrix, residual matrix.}

We extend the function $\w:V\times V\to\bbR_{\geq0}$ to finite paths by setting $\w(u)=\w(x_0,x_1)\times\cdots\times\w(x_{k-1},x_k)$ if $u=(x_0,\ldots,x_k)$, with $\w(u)=1$ if $\length u=0$, and we define for every integer $k\geq0$ and for every pair $(x,y)\in V\times V$:
 \begin{align*}
   Z_{x,y}(k)&=\sum_{u\in O_{x,y}(k)}\w(u),
               &Z_x(k)&=\sum_{u\in O_x(k)}\w(u)=\sum_{y\in V}Z_{x,y}(k)
 \end{align*}
The quantities $Z_{x,y}(k)$ are related to the powers of the adjacency matrix $F$ through the well known formulas:\quad$Z_{x,y}(k)=\evm{ F^k}{x,y}$.

\begin{definition}
  \label{def:1}
  The \emph{growth matrix} $\H(z)$ of a weighted digraph with adjacency matrix $F$ is the power series with matrix coefficients defined by:
  \begin{gather}
    \label{eq:6}
    \H(z)=\sum_{k\geq0}F^kz^k
  \end{gather}
\end{definition}

Just as the adjacency matrix~$F$, the growth matrix $\H(z)$ is defined up to a simultaneous permutation of its lines and columns. The element $\evm{\H(z)}{x,y}$ is itself a generating function:
\begin{gather*}
\evm{\H(z)}{x,y}=\sum_{k\geq0}\evm{F^k}{x,y}\,z^k=\sum_{k\geq0}Z_{x,y}(k)z^k
\end{gather*}

Let $(x,y)\in V\times V$ and let $f(z)=\evm{\H(z)}{x,y}=\sum_ka_kz^k$, say of radius of convergence~$\gamma^{-1}$. It is well known \cite[Ch.V]{flajolet09} that $f(z)$ is a rational series; hence there exists polynomials $A$ and $A_1,\ldots,A_p$ and  complex numbers $\gamma_1,\ldots,\gamma_p$ distinct from $\gamma$ and with $|\gamma_j|\leq|\gamma|$, such that $a_k=\gamma^{k}A(k)+\sum_j\gamma_j^kA_j(k)$ for $k$ large enough. The degree of~$A$, say~$d$, is called the \emph{subexponential degree} of~$f(z)$, and is characterized by the property:
\begin{gather}
  \label{eq:8}
  \lim_{s\to\gamma^{-1}}(1-\gamma s)^{d+1}f(s)=t\in(0,+\infty)
\end{gather}
where the limit is taken for $s\in(0,\gamma^{-1})$; and then the leading coefficient of $A$ is~$\frac 1{d!}t$. When considering the whole matrix~$\H(z)$, we obtain the following result.

\begin{proposition}
  \label{prop:1}
Given a nonnegative matrix~$F$ of spectral radius $\rho>0$ and the matrix power series $\H(z)$ defined as in~\eqref{eq:6}, there exists a unique triple $(\lambda,h,\Theta)$ where $\lambda$ is a positive real, $h$~is a positive integer and $\Theta$ is a square nonnegative matrix, non zero and of the same size as~$F$, such that:
\begin{gather}
  \label{eq:7}
  \lim_{s\to\lambda}(1-\lambda^{-1} s)^h\,\H(s)=\Theta,\qquad s\in(0,\lambda)
\end{gather}

Furthermore, $\lambda=\rho^{-1}$, and $\rho^{-1}$ is the minimal radius of convergence of all the generating functions $\evm{\H(z)}{x,y}$. For each pair $(x,y)$ such that\/ $\evm{\H(z)}{x,y}$ has $\rho^{-1}$ as radius of convergence, then the subexponential degree of\/ $\evm{\H(z)}{x,y}$ is at most~$h-1$, and is equal to $h-1$ if and only if\/ $\evm\Theta{x,y}>0$.
\end{proposition}

\begin{definition}
  \label{def:2}
  Given a weighted digraph $W$ of positive spectral radius and of growth matrix\/~$\H(z)$, the square nonnegative and non zero matrix $\Theta$ as in~\eqref{eq:7} is the \emph{residual matrix} of\/~$W$. The nonnegative integer~$h-1$, where $h>0$ is as in~\eqref{eq:7}, is the \emph{subexponential degree} of\/~$W$.
\end{definition}

\paragraph{\textbullet\bfseries Boltzmann and uniform distributions, weak convergence.}
This paragraph is needed for the reading of Section~\ref{sec:conv-distr}.

\begin{definition}
  \label{def:3}
  Let\/  $W=(V,\w)$ be a weighted digraph. For each vertex $x$ such that $\gamma(x)>0$, we define the family of\/ \emph{Boltzmann distributions} $(\theta_{x,s})_{0<s<\gamma(x)^{-1}}$ and the sequence of \emph{uniform distributions} $(\mu_{x,k})_{k\geq0}$ as follows:
  \begin{gather*}
    \theta_{x,s}=\frac1{G_x(s)}\sum_{u\in O_x}\w(u)s^{\length u}\delta_{\{u\}}\qquad \text{with\quad}G_x(s)=\sum_{u\in O_x}\w(u)s^{\length u}
    \\
     \mu_{x,k}=\frac{1}{Z_x(k)}
\sum_{u\in O_x(k)}\w(u)\delta_{\{u\}}\qquad\text{with\quad}Z_x(k)=\sum_{u\in O_x(k)}\w(u)
  \end{gather*}
where $\delta_{\{u\}}$ is the Dirac distribution concentrated on~$u$.
\end{definition}

The radius of convergence of $\sum_{u\in O_x}\w(u)z^{\length u}$ is $\gamma(x)^{-1}$, hence both $G_x(s)$ and $\theta_{x,s}$ are well defined for $s<\gamma(x)^{-1}$. The assumption $\gamma(x)>0$ is equivalent to the existence of paths in $O_x$ of arbitrary length, hence $Z_x(k)>0$ and $\mu_{x,k}$ is well defined for all $k\geq0$.

All of these  are discrete probability distributions on the set $O_x$ of finite paths starting from~$x$. We can also see them as probability distributions on the set~$\Obar_x$, including infinite paths starting from~$x$. The latter being compact in the product topology, it becomes relevant to look for their weak limits, either for $s\to\gamma(x)^{-1}$ for the Boltzmann distributions, or for $k\to\infty$ for the uniform distributions. Our basic analytical result in this regard will be the following elementary result (for the background on weak convergence, see for instance~\cite{billingsley99}).

\begin{lemma}
  \label{lem:1}
  Let\/ $W$ be a weighted digraph and let $x$ be a vertex.
  \begin{enumerate}
  \item\label{item:5} Assume that, for each integer $k\geq0$, $\mu_k$~is a probability distribution on~$O_x(k)$ and that, for every $u\in O_x$, the following limit exists: $t_u=\lim_{k\to\infty}\mu_k(\up^k u)$.
    Then the sequence $(\mu_k)_{k\geq0}$ converges weakly toward a probability measure $\mu$ on the space $\Omega_x$ of infinite paths starting from~$x$, and $\mu$ is entirely characterized by:\quad
    $  \forall u\in O_x\quad \mu(\up u)=t_u$.
\item\label{item:6} Assume that $(\mu_s)_{0<s<r}$ is a family of probability distributions on~$O_x$ such that, for each $u\in O_x$:\quad
  \begin{inparaenum}[1:]
    \item the following limit exists: $t_u=\lim_{s\to r}\mu_s(\Up u)$; and
    \item $\mu_s(\{u\})\xrightarrow{s\to r}0$.
  \end{inparaenum}
Then the family $(\mu_s)_{0<s<r}$ converges weakly toward a probability measure $\mu$ on~$\Omega_x$ as $s\longrightarrow r$, and $\mu$ is entirely characterized by:\quad     $  \forall u\in O_x\quad \mu(\up u)=t_u$.
  \end{enumerate}
\end{lemma}

\paragraph{\textbullet\bfseries Cocycles and cocycle measures.}
This paragraph is needed for reading Section~\ref{sec:conv-distr} and the end of Section~\ref{sec:umbrella-digraphs}.

Let $W=(V,\w)$ be a weighted digraph and let $\Tos=\{(x,y)\in V\times V\tq x\To y\}$, a subset of $V\times V$ which is of course distinct of $E$ in general. A \emph{cocycle} is a real valued and nonnegative function $\Gamma:\Tos\to\bbR_{\geq0}$ satisfying:
\begin{gather}
  \label{eq:9}
  \forall (x,y,z)\in V\times V\times V\quad (x\To y\To z)\implies \Gamma(x,z)=\Gamma(x,y)\Gamma(y,z)
\end{gather}

Motivated by the form of the limit of the uniform distributions found in~\cite{parry64} for primitive matrices, we introduce below the notion of cocycle measure.

\begin{definition}
  \label{def:4}
  A \emph{cocycle measure} on a weighted digraph\/ $W=(V,\w)$ is a family $\mu=(\mu_x)_{x\in V}$ such that each $\mu_x$ is a probability measure on the space~$\Omega_x$, and such that for some real $\rho>0$ and some cocycle~$\Gamma$, one has:
  \begin{gather}
    \label{eq:10}
    \forall x\in V\quad\forall u\in O_x\quad\mu_x(\up u)=\rho^{-\length u}\w(u)\Gamma(x,\final u)
  \end{gather}
If\/ $\Gamma$ is positive on~$\Tos$, we say that $\mu$ is a \emph{complete cocycle measure}.
\end{definition}

Let $(X_i)_{i\geq0}$ denote the sequence of canonical projections $X_i:\Omega_x\to V$. Assume that $\mu=(\mu_x)_{x\in V}$ is a cocycle measure as in~(\ref{eq:10}), and let $x\in V$. Then, under~$\mu_x$,  $(X_i)_{i\geq0}$~is a Markov chain with initial distribution~$\delta_{\{x\}}$ and with a transition kernel of the following special form, that we call a \emph{cocycle kernel}:
  \begin{gather}
    \label{eq:11}
    q(x,y)=\rho^{-1}\w(x,y)\Gamma(x,y)
  \end{gather}
This correspondence between cocycle measures and cocycle kernels is one-to-one.

\section{Umbrella digraphs}
\label{sec:umbrella-digraphs}

This section is devoted to the study of umbrella digraphs, for which several results from Perron-Frobenius theory for irreducible matrices can be transposed.

\begin{definition}
\label{def:6}
  A weighted digraph\/ $W$ of positive spectral radius is:
  \begin{inparaenum}[\bfseries 1:]
    \item an \emph{umbrella (weighted) digraph} if the basic access classes of\/ $W$ coincide with its final access classes;
    \item an \emph{augmented umbrella (weighted) digraph} if no two distinct basic classes of\/ $W$ have access to each other.
  \end{inparaenum}
\end{definition}


\begin{example}
The digraph of simple elements, but the unit element, of a braid monoid on $n\geq3$ strands is a digraph with two access classes~\cite{abbes19,charney95}. The first class, say~$C_1$, contains only the $\Delta$ element and has access to the second class, say~$C_2$, which contains all the other simple elements. The spectral radii are $\rho(C_1)=1$ and $\rho(C_2)>1$, hence $C_2$ is the unique final and basic class: the digraph is an umbrella digraph.
\end{example}

\begin{example}
  \label{exm:1}
  Let $W=(V,\w)$ be a strongly connected digraph with adjacency matrix $F\neq[0]$. Then $W$ is an umbrella digraph with a unique access class. For $p\geq1$, let $W^{(p)}=(V,\w^{(p)})$ be the digraph with same vertices as~$W$ and with the weight function corresponding the $p^{\text{th}}$~power~$F^p$. Then $W^{(p)}$ is irreducible for all $p\geq1$ if and only if $F$ is primitive, \ie, if $W$ is aperiodic. If $W$ is periodic of period~$d$, then $W^{(d)}$ is an umbrella digraph of spectral radius~$\rho(W)^d$. The digraph $W^{(d)}$ has $d$ basic classes which correspond to the periodic classes of the vertices (see~\cite{seneta81}). All the $d$ access classes of $W^{(d)}$ are both final and basic.
\end{example}

Under an aperiodicity assumption, the adjacency matrix of an augmented umbrella digraph has a decomposition that extends the well known decomposition of primitive matrices recalled in~(\ref{eq:3}).

\begin{theorem}[umbrella and augmented umbrella digraphs with aperiodicity]
  \label{thr:2}
  Let $F$ be the adjacency matrix of an augmented umbrella digraph of spectral radius $\rho>0$ and with $p$ basic classes.  We assume that all the basic classes of\/ $W$ are aperiodic. Then:
  \begin{enumerate}
  \item\label{item:1}
    There is a computable projector matrix of rank $p$ and a matrix $R$ satisfying:
    \begin{align}
      \label{eq:14}
      F&=\rho(\Pi+R)&\Pi\cdot R=R\cdot\Pi&=0&\rho(R)&<1
    \end{align}
    and the following convergence holds:\quad $\lim_{n\to\infty}(\rho^{-1}F)^n=\Pi$.
  \item\label{item:8} There are two computable families of nonnegative line and column vectors $(\ell_i)_{1\leq i\leq p}$ and $(r_i)_{1\leq i\leq p}$ such that:
    \begin{align}
      \label{eq:17}
      \Pi&=\sum_{i=1}^p r_i\cdot\ell_i&\forall i,j\quad \ell_i\cdot r_j&=\delta_i^j
    \end{align}
The family $(r_i)_{1\leq i\leq p}$ is a basis of the space of right $\rho$-eigenvectors of~$F$, and $(\ell_i)_{1\leq i\leq p}$ is a basis of the space of left $\rho$-eigenvectors of~$F$.
\item\label{item:3} Only for umbrella digraphs: for every family $(\alpha_i)_{1\leq i\leq p}$ of reals, the right $\rho$-eigenvector $r=\sum_i\alpha_i r_i$ is positive if and only if $\alpha_i>0$ for all~$i$.
\item\label{item:7} The subexponential degree of\/ $W$ is~$0$. The residual matrix of\/ $W$ coincides with~$\Pi$, \ie: $\lim_{s\to\rho^{-1}}(1-\rho s)\H(s)=\Pi$, where $\H(z)$ is the growth matrix of\/~$W$. Furthermore, $\evm\Pi{x,y}>0$ if and only if there is a basic class $B$ such that $x\To B\To y$.
  \end{enumerate}
\end{theorem}

\begin{remark}
If $W$ is an umbrella digraph, then the condition for $\evm\Pi{x,y}>0$ in point~\ref{item:7} above is equivalent to: $x\To y$ and $y$ belongs to some basic class---the same remark applies to next theorem.
\end{remark}

In the following result, the aperiodicity assumption is dropped. We use the notation $W^{(d)}$ introduced in Example~\ref{exm:1} above.

\begin{theorem}[umbrella and augmented umbrella digraph]
  \label{thr:3}
  Let $F$ be the adjacency matrix of an umbrella (resp., augmented umbrella) digraph\/ $W$ of spectral radius~$\rho>0$.  Let $C_1,\ldots,C_p$ be the basic access classes of\/~$W$, say of periods $d_1,\ldots,d_p$. Let $C_{\nu,0},\ldots,C_{\nu,d_\nu-1}$ be the periodic classes of~$C_\nu$, and let $q=\sum_\nu d_\nu$. Let also $d$ be a common multiple of $d_1,\ldots,d_p$.
  \begin{enumerate}
  \item\label{item:9} $W^{(d)}$ is an umbrella (resp., augmented umbrella) weighted digraph of spectral radius~$\rho^d$ and with $\{C_{\nu,j}\tq 1\leq\nu\leq p,\ 0\leq j<d_\nu\}$ as basic classes, which are all aperiodic. There is a computable projector $\Pi_d$ of rank $q$ and a matrix $R_d$ such that:
    \begin{align}
      \label{eq:18}
      F^d&=\rho^d(\Pi_d+R_d)&\Pi_d\cdot R_d&=R_d\cdot\Pi_d=0&\rho(R_d)&<1
    \end{align}
  \item\label{item:4} There is a basis $(\ell_i)_{1\leq i\leq p}$ of nonnegative left $\rho$-eigenvectors of\/~$F$, and a basis $(r_j)_{1\leq j\leq p}$ of nonnegative right $\rho$-eigenvectors of\/~$F$ satisfying $\ell_i\cdot r_j=\delta_i^j$ for all~$i,j$. 
  \item\label{item:10} Only for umbrella digraphs: for every family $(\alpha_i)_{1\leq i\leq p}$ of reals, the $\rho$-eigenvector $r=\sum_i\alpha_ir_i$ is positive if and only if $\alpha_i>0$ for all~$i$.
  \item\label{item:2} The subexponential degree of\/ $W$ is~$0$. The residual matrix of\/~$W$ is given by
    \begin{gather}
      \label{eq:19}
      \Theta=\frac1d\Bigl(\,\sum_{i=0}^{d-1}\rho^{-i}F^i\Bigr)\cdot\Pi_d
    \end{gather}
    and satisfies: $\evm\Theta{x,y}>0$ if and only if there is a basic class $B$ such that $x\To B\To y$.
  \end{enumerate}
\end{theorem}

In the particular case of a periodic and strongly connected digraph, the following result shows that some simplifications occur in~(\ref{eq:19}); it is the basis of the convergence result for strongly connected digraphs, Th.~\ref{thr:9} in Sect.~\ref{sec:itsh-unif-distr}.

\begin{corollary}[residual matrix of a strongly connected digraph]
  \label{cor:2}
  Let $F$ be the adjacency matrix of a strongly connected digraph of positive spectral radius. Then the subexponential degree of\/ $W$ is~$0$. Let $(\ell,r)$ be a pair of left and right Perron eigenvectors of\/~$F$ such that $\ell\cdot r=1$. Then the residual matrix of\/ $W$ is $\Theta=\Pi$, the positive rank $1$ projector given by $\Pi=r\cdot\ell$.
\end{corollary}

Finally, umbrella weighted digraphs can be characterized by the existence of \emph{complete} cocycle measures.

\begin{theorem}[existence of a complete cocycle measure]
  \label{thr:6}
A digraph\/ $W$ is an umbrella digraph if and only if there exists a complete cocycle measure on~$W$.

If\/ $W$ is an umbrella digraph and is accessible from a vertex~$x$, then the complete cocycle measures on\/ $W$ are parameterized by an open simplex of dimension~$p-1$, where $p$ is the number of basic classes of\/~$W$.
\end{theorem}

In particular there exists a unique complete cocycle measure if $W$ is strongly connected. For $F$ the adjacency matrix with $\rho=\rho(F)$, and for $r$ a Perron right eigenvector of~$F$, the cocycle kernel associated to this unique complete cocycle measure is given by:
\begin{align*}
  q(x,y)&=\rho^{-1}\w(x,y)\Gamma(x,y)&\Gamma(x,y)&=\frac{\evm ry}{\evm rx}
\end{align*}

The following well known result~\cite{rothblum14,berman94} can be seen as a consequence of Th.~\ref{thr:6}.

\begin{corollary}[existence of a positive Perron eigenvector]
  \label{cor:1}
  A weighted digraph $W$ of spectral radius $\rho>0$ is an umbrella weighted digraph if and only if its adjacency matrix has a positive $\rho$-eigenvector.
\end{corollary}

\section{Computing the residual matrix}
\label{sec:comp-asympt-matr}

This section is devoted to the recursive computation of the residual matrix of a general weighted digraph.



\smallskip
Recall that a subset $A$ of a poset $(B,\leq)$ is:
\begin{inparaenum}[\bfseries 1:]
  \item \emph{final} if:\quad$\forall (a,b)\in A\times B\quad a\leq b\implies b\in A$; and
  \item \emph{initial} if:\quad$\forall (a,b)\in A\times B\quad b\leq a\implies b\in A$.
\end{inparaenum}

\subsubsection{\textbullet\itshape Theoretical results.}
\label{sec:theoretical-results}

The following elementary result will be instrumental.

\begin{lemma}[recursive form of the growth matrix]
  \label{lem:2}
  Let\/ $W=(V,\w)$ be a weighted digraph with adjacency matrix~$F$. Assume that\/ $(S,T)$ is a partition of\/ $V$ and that $T$ is final in~$V$. Let\/ $\H_S(z)$ and\/ $\H_T(z)$ be the growth matrices of\/~$S$ and\/~$T$, sub-weighted digraphs of\/~$W$. Then:
  \begin{align}
    \label{eq:4}
    \H(z)&=\begin{pmatrix}
      \H_S(z)&\Y(z)\\
      0&\H_T(z)
    \end{pmatrix}
  &\Y(z)&=z\,\H_S(z)\cdot X\cdot\H_T(z)
\end{align}
where $X$ is the rectangular block submatrix of\/ $F$ corresponding to~$S\times T$.
\end{lemma}

To determine the subexponentiel degree of a general weighted digraph, the notion of height that we introduce below and which is well known in nonnegative matrix theory~\cite{rothblum14,berman94}, plays a key role.

\begin{definition}
  \label{def:5}
  Let\/ $W$ be a weighted digraph, and let\/ $(\D,\To)$ be the poset of its access classes. A \emph{dominant chain} is a chain of basic classes in $(\D,\To)$ and of maximal length. The \emph{height} of\/ $W$ is the length of the dominant chains.
\end{definition}


\begin{remark}
  \label{rem:1}
  The augmented umbrella digraphs from Def.~\ref{def:6} are the weighted digraphs of positive spectral radius and of height~$1$.
\end{remark}

\begin{theorem}[subexponential degree of a digraph]
  \label{thr:1}
  Let\/ $W=(V,\w)$ be a weighted digraph of positive spectral radius and of height~$h$.
  \begin{enumerate}
  \item The subexponential degree of\/ $W$ is $h-1$.
  \item Let\/ $\Theta$ be the residual matrix of\/~$W$. Then $\evm\Theta{x,y}>0$ if and only if there exists a dominant chain $(L_1,\ldots,L_p)$ such that $x\To L_1$ and $L_p\To y$.
  \end{enumerate}
\end{theorem}

From the above result derives the following well knwown result originally proved by Rothblum~\cite{rothblum14}.

\begin{corollary}
  \label{cor:3}
  Let $F$ be the adjacency of a weighted digraph of spectral radius $\rho$ and of height~$h$. Then the dimension of the generalized eigenspace of $F$ associated to $\rho$ is~$h$.
\end{corollary}

\subsubsection{\textbullet\itshape Recursive computing.}
\label{sec:recursive-computing}

We aim at recursively computing the residual matrix of a weighted digraph~$W$ of spectral radius $\rho>0$. We first observe two facts.

\textsc{Fact~1.}\quad The growth matrix $\H(z)$ can be recursively computed by starting from its lower-right corner and extending the blocks already computed. In details, let $(D_1,\ldots,D_p)$ be an enumeration of the access classes of~$W$. Put $V_i=D_1\cup\ldots\cup D_i$ for $i=0,\ldots,p$, and assume that the enumeration has been chosen such that $V_i$ is final in $V$ for each~$i$. Then, with the obvious notations $\H_{D_i}(z)$ and~$\H_{V_i}(z)$, one has for each~$i>0$:
\begin{align}
\label{eq:13}
  \H_{V_i}(z)&=\begin{pmatrix}
      \H_{D_i}(z)&\Y_i(z)\\
      0&\H_{V_{i-1}}(z)
    \end{pmatrix}
  &\Y_i(z)&=z\,\H_{D_i}(z)\cdot X_i\cdot\H_{V_{i-1}}(z)
\end{align}
where $X_i$ is the rectangular block submatrix of $F$ corresponding to $D_i\times{V_{i-1}}$. This results from Lemma~\ref{lem:2} applied with the partition $(D_i,V_{i-1})$ of~$V_i$.\eof

\textsc{Fact~2.}\quad If the height of $W$ is~$1$, then the residual matrix of $W$ is directly computable. This results from Th.~\ref{thr:3}, \emph{via} Remark~\ref{rem:1} above.\eof

\smallskip
Let $r_i=\rho(D_i)$ for $0<i\leq p$. To simplify the exposition, we assume that $r_i>0$ for all~$i$, and we consider $\Pi_i$ the  residual matrix of~$D_i$, which is a computable rank~$1$ projector (see Corollary~\ref{cor:2}).

\smallskip
\textsc{Initialization.}\quad Denoting $\rho_i=\rho(V_i)$, we start the induction with the first $i_0$ such that $h_{i_0}\geq2$ and $\rho_{i_0}=\rho$; in particular, $\rho(V_{i_0-1})=\rho$. Indeed, as long as the height of $V_i$ is~$1$, its residual matrix can be directly computed, as we observed in Fact~2 above.

Let  $\Theta_i$ be the residual matrix of~$V_i$, and let $h_i$ denote the height of~$V_i$. Since $\rho(V_i)=\rho$ for all $i\geq i_0$, we have:
\begin{gather}
  \label{eq:12}
\forall i\geq i_0\quad \Theta_i=\lim_{s\to\rho^{-1}}(1-\rho s)^{h_i}\,\H_{V_i}(s)
\end{gather}

\smallskip
\textsc{Induction step.} We assume that the residual matrix $\Theta_{i-1}$ of $V_{i-1}$ has been computed, and in certain cases, we also need to call for the computation of the residual matrix of some sub-digraph of~$V_{i-1}$. We show how to compute~$\Theta_i$.
\begin{enumerate}
\item \emph{Case where $h_i=h_{i-1}$.}
  \begin{enumerate}
  \item \emph{If $\rho(D_i)<\rho$.}\quad Referring to~\eqref{eq:13} and~\eqref{eq:12}, and since $h_i\geq2$, we have:
    \begin{gather*}
      (1-\rho s)^{h_i}\,   \Y_i(s)=s\, \H_{D_i}(s)\cdot X_i\cdot\bigl((1-\rho s)^{h_{i-1}}\,\H_{V_{i-1}}(s)\bigr)\\
\begin{aligned}
\text{and thus\quad}
      \Theta_i&=
      \begin{pmatrix}
        0&\ & A\\
        0&\  &\Theta_{i-1}
      \end{pmatrix}
          &\text{with\quad}A&=\rho^{-1}\,\H_{D_i}(\rho^{-1})\cdot X_i\cdot\Theta_{i-1}
        \end{aligned}
      \end{gather*}
We have  $\H_{D_i}(\rho^{-1})=(\mathrm{Id}-\rho^{-1}F_i)^{-1}$, where $F_i$ is the block sub-matrix of $F$ corresponding to $D_i\times D_i$; so $\Theta_i$ can be computed.
    \item \emph{If $\rho(D_i)=\rho$}.\quad Consider the partition $(\Ttilde,\Tbar)$ of~$V_{i-1}$ with
    \begin{gather}
      \label{eq:16}
      \Ttilde    =\{y\in V_{i-1}\tq\exists x\in V_{i-1}\quad (y\To x)\land( D_i\To x)\}
    \end{gather}
which is initial in~$V_{i-1}$. Enumerating the vertices of $V_i$ as those of~$D_i$, then those of~$\Ttilde$ and then those of~$\Tbar$ in this order, we have:
      \begin{align}
\label{eq:15}
        \H_{V_i}(z)&=\begin{pmatrix}
          \H_{D_i}(z)&\A(z)&0\\
          0&\H_{\Ttilde}(z)&\B(z)\\
          0&0&\H_{\Tbar}(z)
        \end{pmatrix}&
                         \A(z)&=z\,\H_{D_i}(z)\cdot A\cdot\H_{\Ttilde}(z)
      \end{align}
      where $A$ and $B$ are the block submatrices of $F$ corresponding to $D_i\times\Ttilde$ and to~$\Ttilde\times\Tbar$. Indeed, vertices of~$\Tbar$ are not accessible from~$D_i$, whence the zero block in the upper right corner of~$\H_{V_i}(z)$.

      On the one hand, the down right corner of $\H_{V_i}(z)$ formed by the four blocks in~(\ref{eq:15}) is nothing but $\H_{V_{i-1}}(z)$. On the other hand, we have:
      \begin{align*}
        (1-\rho s)^{h_i}\,\A(s)&=s\,\bigl((1-\rho s)\H_{D_i}(s)\bigr)\cdot A\cdot\bigl((1-\rho s)^{h_i-1}\,\H_{\Ttilde}(s)\bigr)\\
                               &\xrightarrow{s\to\rho^{-1}}\Y
=\begin{cases}
                                 0,&\text{if the height of $\Ttilde$ is $< h_i-1$}\\
                                 \rho^{-1}\,\Pi_i\cdot A\cdot \Ytilde&\text{if the height of $\Ttilde$ is $h_i-1$}
                               \end{cases}
      \end{align*}
      where $\Ytilde$ is the residual matrix of~$\Ttilde$. Hence $\Theta_i=\left(\begin{smallmatrix}
          0&\Y\quad 0\\0&\Theta_{i-1}
          \end{smallmatrix}\right)$ is computable. 
    \end{enumerate}
  \item \emph{Case where $h_i=h_{i-1}+1$.} According to~\eqref{eq:13}, we compute as follows:
    \begin{align*}
      (1-\rho s)^{h_i}\,\Y(s)&=s\,\bigl((1-\rho s)\H_{D_i}(s)\bigr)\cdot X_i
                               \cdot\bigl((1-\rho s)^{h_{i-1}}\,\H_{V_{i-1}}(s)\bigr)\\
      &\xrightarrow{s\to\rho^{-1}}\Y=\rho^{-1}\,\Pi_i\cdot X_i\cdot\Theta_{i-1}
    \end{align*}
Hence $\Theta_i=\left(
  \begin{smallmatrix}
    0&\Y\\0&0
  \end{smallmatrix}
\right)$ is computable.
\end{enumerate}

\section{Convergence of distributions}
\label{sec:conv-distr}

This section is devoted to the convergence of distributions on paths of a weighted digraph, which was our main goal from the beginning. We focus on the family $(\theta_{x,s})_{0<s<\rho^{-1}}$ of Boltzmann distributions on the one hand, and on the sequence $(\mu_k)_{k\geq0}$ of uniform distributions on the other hand (see Def.~\ref{def:3}).

\subsubsection{\textbullet\itshape Boltzmann distributions.}
\label{sec:itsh-boltzm-distr}

Let $W=(V,\w)$ be a weighted digraph of spectral radius $\rho>0$. In view of Lemma~\ref{lem:1}, point~\ref{item:6}, we fix a pair $(x,u)$ where $x\in V$ and $u\in O_x$, and we study the quantity $\theta_{x,s}(\Up u)$, aiming at its convergence for $s\to\rho^{-1}$. For $s<\rho^{-1}$, we have:
\begin{align}
  \label{eq:20}
  \theta_{x,s}(\Up x)&=\frac1{G_x(s)}\sum_{v\in O_x\tq u\leq v}\w(v)s^{\length v}
\end{align}
Every $v\in O_x$ with $u\leq v$ writes in a unique way as the concatenation $v=u\cdot v'$ for some $v'\in O_{\final u}$, and then $\w(v)=\w(u)\w(v')$. Hence (\ref{eq:20}) writes as:
\begin{align}
  \label{eq:23}
  \theta_{x,s}(\Up x)&=\frac{\w(u)s^{\length u}}{G_x(s)}\sum_{v\in O_{\final u}}\w(v)s^{\length v}
                       =\w(u)s^{\length u}\frac{G_{\final u}(s)}{G_x(s)}
\end{align}
Now for every vertex $y$ and for every real $s\in(0,\rho^{-1})$, $G_y(s)$~is related to the growth matrix of $W$ \emph{via}:
\begin{align}
  \label{eq:22}
  G_y(s)&=\sum_{t\in V}\evm{\H(s)}{y,t}=\evm{\H(s)\cdot 1}{y}
\end{align}
where $1$ denotes the column vector filled with~$1$s.

Let $h(x)$ be the height of $V(x)=\{y\in V\tq x\To y\}$, and let $\Theta$ be the residual matrix of~$V(x)$. It follows from Th.~\ref{thr:1} that an equivalent of $\evm{\H(s)\cdot 1}x$ when $s\to\rho^{-1}$ is:
\begin{gather}
  \label{eq:21}
  \evm{\H(s)\cdot 1}x\sim_{s\to\rho^{-1}} \bigl(h(x)\bigr)!\,(1-\rho s)^{-h(x)}\, \evm{\Theta\cdot 1}x
\end{gather}
since $\evm{\Theta\cdot 1}x>0$, as a sum of terms of which at least one is positive according to Th.~\ref{thr:1}. Putting $y=\final u$, we note that $y\in V(x)$ and that two cases may occur.
\begin{inparaenum}[\bfseries 1:]
  \item If there exists a dominant chain $(L_1,\ldots,L_{h(x)})$ of $V(x)$ such that $y\To L_1$, then an equivalent for $\evm{\H(s)\cdot1}y$ analogous to~(\ref{eq:21}) and with the same exponent $h(x)$ holds. From~\eqref{eq:23},~\eqref{eq:22} and~\eqref{eq:21} we derive thus:
    \begin{gather}
      \label{eq:24}
      \lim_{s\to\rho^{-1}}\theta_{x,s}(\Up u)=\rho^{-\length u}\w(u)\frac{\evm{\Theta\cdot 1}{\final u}}
      {\evm{\Theta\cdot 1}x}
    \end{gather}
    which is a positive number.
    
But, 
\item If not, it means that $\evm{\H(s)\cdot 1}y$ has either an exponential growth rate less that~$\rho$, or an exponential growth rate equal to $\rho$ but a subexponential degree less than~$h$. In both situations, the ratio in~(\ref{eq:23}) goes to zero as $s\to\rho^{-1}$.
\end{inparaenum}

\smallskip
The above discussion motivates the following definition, where we recall that $\gamma(x)=\rho\bigl(V(x)\bigr)$.

\begin{definition}
  \label{def:7}
  Let $W=(V,\w)$ be a weighted digraph and let $x$ be a vertex such that $\gamma(x)>0$. The \emph{umbrella digraph spanned by~$x$} is the sub-digraph of\/ $V(x)$, the vertices of which are the vertices $y\in V(x)$ such that $y\To L_1$ for some dominant chain $(L_1,\ldots,L_{h(x)})$ of\/~$V(x)$. We denote it by~$\U(x)$.
\end{definition}

The digraph $\U(x)$ thus defined is indeed an umbrella digraph, of spectral radius~$\gamma(x)$. Putting together the result of the above discussion and Lemma~\ref{lem:1}, we obtain the following convergence result.

\begin{theorem}[convergence of Boltzmann distributions]
  \label{thr:7}
  Let\/ $W$ be a weighted digraph of positive spectral radius~$\rho$, with residual matrix~$\Theta$, and let $x_0$ be a vertex such that all vertices of\/ $W$ are accessible from~$x_0$. Then the Boltzmann distributions\/ $\theta_{x_0,s}$ converge weakly when $s\to\rho^{-1}$ toward the complete cocycle measure $\theta$ on $\U(x_0)$ which cocycle transition kernel is given on $\U(x_0)$ by:
  \begin{gather}
    \label{eq:25}
    q(x,y)=\rho^{-1}\w(x,y)\Gamma(x,y)\qquad\Gamma(x,y)=\frac{\evm{\Theta\cdot 1}y}{\evm{\Theta\cdot1}x}
  \end{gather}
\end{theorem}

\begin{remark}
  \label{rem:2}
  In the above result, the condition that $V$ should be accessible from $x_0$ is not a severe restriction. In general, the theorem applies to the sub-digraph~$V(x_0)$, which vertices are those accessible from~$x_0$ and of spectral radius~$\gamma(x_0)$. If $\gamma(x_0)=0$, the theorem does not apply; and indeed, $\Omega_{x_0}=\emptyset$ in this case.
\end{remark}

\begin{remark}
  \label{rem:3}
  The limit distributions $(\theta_x)_{x\in V}$ are also obtained as the limits of the Boltzmann distributions relative to the sub-digraph $\U(x_0)$ itself. Hence if $\widetilde\Theta$ is the residual matrix of~$\U(x_0)$, the cocycle kernel $q(x,y)$ in~\eqref{eq:25} can also be obtained as $q(x,y)=\rho^{-1}\w(x,y)\widetilde\Gamma(x,y)$ with $\widetilde\Gamma(x,y)=\evm{\widetilde\Theta\cdot1}{y}/\evm{\widetilde\Theta\cdot1}y$. And $\widetilde\Theta$ is directly computed according to Th.~\ref{thr:3} since $\U(x_0)$ is an umbrella digraph.
\end{remark}

\subsubsection{\textbullet\itshape Uniform distributions.}
\label{sec:itsh-unif-distr}

We now aim at using point~\ref{item:5} of Lemma~\ref{lem:1} in order to derive the weak convergence of the sequence $(\mu_{x,k})_{k\geq0}$ of uniform distributions on paths. We fix a vertex $x$ and a path $u\in O_x$, and we consider the quantity $\mu_{x,k}(\up^k u)$ for $k\geq\length u$. The same change of variable that we used above yields the following expression:
\begin{gather}
  \label{eq:26}
  \mu_x(\up^k u)=\frac1{Z_x(k)}\Bigl(\,\sum_{v\in O_x(k)\tq u\leq v}\w(v)\Bigr)=
\w(u)\frac{Z_{\final u}\bigl(k-\length u\bigr)}{Z_x(k)}
\end{gather}

We are thus brought to discuss the asymptotics of the $k^\text{th}$ coefficients of the generating functions $\evm{\H(z)\cdot 1}x$ and $\evm{\H(z)\cdot 1}{\final u}$, as $k\to\infty$. These are closely related to the values of these generating functions near their singularity~$\rho^{-1}$, which connects to our previous discussion for the Boltzmann distributions. We have thus the following result, which is a sort of Tauberian theorem relatively to Th.~\ref{thr:7}. Remarks~\ref{rem:2} and~\ref{rem:3} above apply to Th.~\ref{thr:8} as they did for Th.~\ref{thr:7}.

\begin{theorem}[convergence of uniform distributions with aperiodicity]
  \label{thr:8}
  Let\/ $W$ be a weighted digraph of positive spectral radius~$\rho$, and let $x_0$ be a vertex such that all vertices of\/ $W$ are accessible from~$x_0$. We assume that all the basic classes of\/ $W$ are aperiodic. Then the sequence $(\mu_{x_0,k})_{k\geq0}$ of uniform distributions converges weakly toward the complete cocycle measure on $\U(x_0)$ with the cocycle transition kernel $q$ described in~\eqref{eq:25}.
\end{theorem}

The aperiodicity condition in the above theorem is sufficient but not necessary, as shown by the following result, consequence of the simple form of the residual matrix for a strongly connected digraph as given in Cor.~\ref{cor:2}.

\begin{theorem}[convergence of uniform distributions for a strongly connected digraph]
  \label{thr:9}
Let\/ $W=(V,\w)$ be a strongly connected digraph with positive spectral radius~$\rho$. Then $\U(x)=V$ for every vertex~$x$. The sequence $(\mu_{x,k})_{k\geq0}$ of uniform distributions converges weakly toward the unique complete cocycle measure on~$W$. The corresponding transition kernel is $q(x,y)=\rho^{-1}\w(x,y)\Gamma(x,y)$ where the cocycle $\Gamma$ is given by $\Gamma(x,y)=\evm ry/\evm rx$ for any Perron eigenvector $r$ of the adjacency matrix of\/~$W$.
\end{theorem}

In general, the convergence of uniform distributions does not hold. In the case of an umbrella weighted digraph, we have the following condition.

\begin{theorem}[convergence of uniform distributions for an umbrella digraph]
  \label{thr:10}
  Let\/ $W=(V,\w)$ be an umbrella weighted digraph of spectral radius~$\rho$ and with adjacency matrix~$F$. Let $d$ be a common multiple of the periods of the basic classes of\/~$W$. Let $\Pi_d$ be the projector defined as in Th.~\ref{thr:3} and define $\beta_{x,i}=\rho^{-i}\evm{\Pi_d\cdot F^i\cdot 1}x$ for every vertex $x$ and every integer $0\leq i<d$. 

    Then $(\mu_{x,k})_{k\geq0}$ has a weak limit if and only if $\beta_{x,i}$ is independent of~$i$; this condition is automatically satisfied for the vertices $x$ belonging to the basic classes of\/~$W$.
\end{theorem}


%

\printbibliography

\newpage
\appendix

\section{Proofs for Section~\ref{sec:umbrella-digraphs}}

\subsection{Proof of Theorem~\ref{thr:2}}

We give the proof only for the case of umbrella digraphs. The case of augmented umbrella digraphs is similar.
We organize the vertices of $W$ as follows:
\begin{itemize}
\item[$T_1$ :] the vertices which have access to some basic class without being in a basic class
\item[$T_2$ :] the vertices in some basic class
\end{itemize}

Accordingly, the adjacency matrix $F$ has the following block decomposition:
 \begin{gather}
   \label{eq:27}
   F=
   \begin{pmatrix}
     A&X_1&\hdotsfor{2}& X_p\\
     0&F_1&0&\dots&0\\
     \vdots&
     &&&\vdots\\
     0&\hdotsfor{2}&0&F_p
   \end{pmatrix}
 \end{gather}
The diagonal blocks $F_1,\ldots,F_p$ are the adjacency matrices of the basic classes, $A$~is a square nonnegative matrices with $\rho(A)<\rho$, and $X_1,\ldots,X_p$ are nonnegative rectangular block matrices. 

For each $i\in\{1,\ldots,p\}$, $F_i$~is irreducible and $\rho(F_i)=\rho$. Let $(\elltilde_i,\rtilde_i)$ be a pair of left and right positive $\rho$-eigenvectors of $F_i$ satisfying $\elltilde_i\cdot \rtilde_i=1$, which exists according to Perron-Frobenius theory for irreducible matrices. We put:
\begin{align}
  \label{eq:28}
  r_i&=\begin{pmatrix}
u_i\\0\\
\rtilde_i\\
0
\end{pmatrix}
               &
                 \ell_i&=
                        \begin{pmatrix}
                          0&\dots\elltilde_i&\dots&0
                        \end{pmatrix}
\end{align}
where $u_i$ is to be chosen so as $r_i$ being a $\rho$-eigenvectors of~$F$. For $r_i$ to be a right $\rho$-eigenvector of~$F$, it is necessary and sufficient that: $A\cdot u_i+X_i\cdot\rtilde_i=\rho u_i$, hence: $(\rho\Id-A)\cdot u_i=X_i\cdot\rtilde_i$. Since $\rho(A)<\rho$, there is a unique solution, given by:
\begin{gather}
  \label{eq:29}
  u_i=(\rho\Id-A)^{-1}\cdot X_i\cdot\rtilde_i=\sum_{k\geq0}\rho^{-(k+1)}A^k\cdot X_i\cdot\rtilde_i
\end{gather}

On the other hand, $\ell_i$ defined in~(\ref{eq:28}) is directly a left $\rho$-eigenvector of~$F$. It is apparent on~(\ref{eq:28}) that: $\ell_i\cdot r_i=\elltilde_i\cdot\rtilde_i=1$, and that $\ell_i\cdot r_j=0$ if $i\neq j$.

The expression on the rightmost members of~(\ref{eq:29}) shows that $u_i$ is a nonnegative vector, since all the vectors and matrices involved are nonnegative. Hence $r_i$ and $\ell_i$ are nonnegative vectors. Furthermore, $F_i$~being irreducible, $\rho$~has algebraic multiplicity $1$ in $F_i$ according to Perron-Frobenius theory for irreducible matrices. Hence $\rho$ has algebraic multiplicity $p$ in~$F$ and  $(r_i)_{1\leq i\leq p}$ is a basis of right $\rho$-eigenvectors of~$F$ and $(\ell_i)_{1\leq i\leq p}$ is a basis of left $\rho$-eigenvectors of~$F$.

We define $\Pi$ and $R$ as follows:
\begin{align}
  \label{eq:30}
  \Pi&=\sum_{i=1}^p r_i\cdot\ell_i&R=\rho^{-1}F-\Pi
\end{align}
The following computations show that $\Pi$ is a projector and that $\Pi\cdot R=R\cdot\Pi=0$:
\begin{align*}
  \Pi^2&=\sum_{1\leq i,j\leq p}r_i\cdot(\ell_i\cdot r_j)\cdot\ell_j=\sum_{i=1}^pr_i\cdot\ell_i=\Pi\quad\text{using $\ell_i\cdot r_j=\delta_i^j$}\\
  \Pi\cdot R&=\rho^{-1}\Pi\cdot F-\Pi^2=\sum_{i=1}^pr_i\cdot(\rho^{-1}\ell_i\cdot F)-\Pi=\sum_{i=1}^pr_i\cdot\ell_i-\Pi=0\\
  R\cdot\Pi&=\rho^{-1}F\cdot\Pi-\Pi^2=\sum_{i=1}^p(\rho^{-1}F\cdot r_i)\cdot\ell_i-\Pi=\sum_{i=1}^pr_i\cdot\ell_i-\Pi=0
\end{align*}
By its definition~(\ref{eq:30}), one sees that $\Pi$ is indeed of rank~$p$.

To show that $\rho(R)<1$, we observe that $\Pi$ and $R$ have the following block decompositions:
\begin{align*}
  \Pi&=
       \begin{pmatrix}
         0&\X&\hdotsfor{1}&\X\\
         0&\Pi_1&\dots&0\\
         \vdots&\ddots&\ddots&\vdots\\
         0&\dots&0&\Pi_p
       \end{pmatrix}
&R&=
    \begin{pmatrix}
      \rho^{-1}A&\X&\hdotsfor{1}&\X\\
         0&R_1&\dots&0\\
         \vdots&\ddots&\ddots&\vdots\\
         0&\dots&0&R_p
    \end{pmatrix}
&
                    \begin{aligned}
R_i&=\rho^{-1}F_i-\Pi_i\\
\Pi_i&=\rtilde_i\cdot\elltilde_i                      
                    \end{aligned}
\end{align*}
where the $\X$s represent rectangular nonnegative block matrices. Hence to prove that $\rho(R)<1$, it is enough to see that $\rho(R_i)<1$ for all~$i$. But this derives from the standard Perron-Frobenius theory for primitive matrices, since $F_i$ is indeed primitive thanks to the aperiodicity assumption in the statement.

As a consequence of $F=\rho^{-1}(\Pi+R)$ on the one hand, and of $\Pi\cdot R=R\cdot\Pi=0$ and $\Pi^2=\Pi$ on the other hand, we have:
\begin{gather}
  \label{eq:32}
\forall k>0\quad (\rho^{-1} F)^k=\Pi+R^k  
\end{gather}
and $R^k\xrightarrow{k\to\infty}0$ since $\rho(R)<1$, hence $(\rho^{-1}F)^k\xrightarrow{k\to\infty}\Pi$. So far we have proved points~\ref{item:1} and~\ref{item:8} of Th.~\ref{thr:2}.

We now prove point~\ref{item:7}. \emph{Residual matrix~$\Theta=\Pi$}.
From~(\ref{eq:32}), we get:
\begin{align}
  \label{eq:37}
  \H(z)&=\Id+\frac{\rho z}{1-\rho z}\Pi+\T(z)&\T(z)&=\sum_{k>0}(\rho R)^kz^k
\end{align}
Since $\rho(R)<1$, there is a matrix norm $\|\cdot\|$ such that $\|R\|<1$, and for this matrix norm:
\begin{gather}
  \notag
  \|\T(z)\|\leq\frac{\|\rho R z\|}{1-\|\rho Rz\|}\\
  \label{eq:36}
s\in(0,\rho^{-1})\qquad \|  (1-\rho s)\T(s)\|\leq (1-\rho s)\|\frac{\|R\|}{1-\|R\|}
\end{gather}
Putting together~(\ref{eq:37}) and~(\ref{eq:36}) yields:
\begin{gather*}
  \lim_{s\to\rho^{-1}}(1-\rho s)\H(s)=\Pi
\end{gather*}
and since $\Pi$ is non zero, that proves that the subexponential degree of $W$ is $0$ and that $\Theta=\Pi$.

\emph{Condition for $\evm\Pi{x,y}>0$}. Consider the following equivalence:
\begin{gather}
  \label{eq:38}
\bigl(\evm\Pi{x,y}>0\bigr)\iff\bigl((x\To y)\land( y\in T_2)\bigr)
\end{gather}
The equivalence is true if $x$ and $y$ belong to~$T_2$, since both terms of~\eqref{eq:38} are true. The equivalence is also true if $y$ belongs to~$T_1$ since both terms of~\eqref{eq:38} are false. It remains only to settle the case where $x\in T_1$ and $y\in T_2$. In this case, and if $x\To y$, then the right term of~\eqref{eq:38} is true. We show that the left term is also true. First, observe that $\evm\Pi{x,y}=\evm{u_i}x\evm{\elltilde_i}y$, where $i$ the index of the basic class to which $y$ belongs to. Since $\elltilde_i>0$, is remains only to show that $\evm{u_i}x>0$. But the assumption $x\To y$ implies that, for some integer~$k>0$: $\evm{F^k}{x,y}>0$, which implies that $\evm{A^k\cdot X_i}{x,y}>0$. By the form~\eqref{eq:29} of~$u_i$, and since $\rtilde_i>0$, we deduce that $\evm{u_i}x>0$, which was to be shown. And if $\neg(x\To y)$, then both terms of~\eqref{eq:38} are false, hence the equivalence~\eqref{eq:38} is true in all cases.

Finally we prove point~\ref{item:3}. Let $(\alpha_i)_{1\leq i\leq p}$ be a collection of reals and let $r=\sum_i\alpha_ir_i$ which is a $\rho$-eigenvector of~$F$. Assuming that $r>0$, we prove that $\alpha_i>0$ for all~$i$. Let $i\in\{1,\ldots,p\}$, and pick $x$ in the correspond basic class. Then $\evm rx=\alpha_i\evm{\rtilde_i}x$ and thus $\alpha_i>0$.

Conversely, assume that $\alpha_i>0$ for all~$i$. For $x\in T_2$, say belonging to the $i^\text{th}$ basic class, we have $\evm rx=\alpha_i \evm{\rtilde_i}x>0$ since $\rtilde_i>0$. For $x\in T_1$, aiming at proving that $\evm rx>0$, let $u=\sum_i\alpha_i u_i$. By construction, $x$~has access to at least one vertex $y\in T_2$, say belonging to the $i^\text{th}$ basic class. Therefore $\evm{A^k\cdot X_i}{x,y}>0$ for some integer $k>0$, and since $\rtilde_i>0$, (\ref{eq:29})~implies that $\evm{u_i}x>0$ and thus $\evm ux>0$. Finally $\evm rx=\evm ux>0$, which completes the proof.

\subsection{Proof of Theorem~\ref{thr:3}}
\label{sec:proof-theorem}

Point~\ref{item:9} derives from the theory of periodic irreducible matrices, see~\cite{seneta81}, and from Th.~\ref{thr:2} applied to~$W^{(d)}$. The construction of $(\ell_i)_{1\leq i\leq p}$ and of $(r_i)_{1\leq i\leq p}$ is the same as in the proof of Th.~\ref{thr:2}; indeed, only the irreducibility of the basic classes was used at this stage, not their aperiodicity. The same applies for the proof of point~\ref{item:10} on the positivity condition of the $\rho$-eigenvector $\sum_i\alpha_i r_i$ if $W$ is an umbrella digraph.

It remains only to prove point~\ref{item:2} by computing the residual matrix of~$W$. Let $\H(z)$ and $\H^{(d)}(z)$ be the growth matrices of $W$ and of~$W^{(d)}$. Then:
\begin{align*}
  \H(z)&=\Bigl(\,\sum_{i=0}^{d-1}z^iF^i\Bigr)\cdot\H(z^d)&
                                                           \lim_{s\to\rho^{-d}}(1-\rho^ds)\,\H^{(d)}(s)&=\Pi_d
\end{align*}
Taking the limit of $(1-\rho s)\H(s)$ yields the result for~$\Theta$, observing the occurrence of the factor~$\frac 1d$.

\subsection{Proof of Corollary~\ref{cor:2}}

Let $F^d=\rho^d(\Pi_d+R)$ be the decomposition of the adjacency matrix $F^d$ of $W^{(d)}$ given by Th.~\ref{thr:2} with $\Pi_d\cdot R=R\cdot\Pi_d=0$, $\Pi_d$~a projector and $\rho(R)<1$.
From the theory of irreducible matrices, we now that, after a simultaneous permutation of the lines and the columns, $F$~and $F^d$ can be written by blocks as follows:
\begin{align*}
    F&=
     \begin{pmatrix}
       0&Q_0&\ldots&0\\
       \vdots&\ddots&\ddots&0\\
       &&0&Q_{d-2}\\
       Q_{d-1}&&\ldots&0
     \end{pmatrix}
&F^d&=
      \begin{pmatrix}
        S_0&\ldots&0\\
        \vdots&\ddots&\vdots\\
        0&\ldots&S_{d-1}
      \end{pmatrix}
\end{align*}
with $S_i=Q_i\cdot\ldots Q_{i-1}$ (a circular product of $d$ matrices), all $S_i$ are primitive of spectral radius~$\rho^d$. Furthermore, if $(\elltilde_0,\rtilde_0)$ is a pair of positive left and right $\rho^d$-eigenvectors of $S_0$ such that $\elltilde_0\cdot\rtilde_0=1$, an analogous pair for $S_i$ is given by $\elltilde_i=\rho^{-i}\elltilde_0\cdot Q_0\cdot\ldots\cdot Q_{i-1}$ and $\rtilde_i=\rho^{-d+i} Q_i\cdot\ldots\cdot Q_{d-1}\cdot \rtilde_0$\,. Additionally, putting $\ell_i=\left(\begin{smallmatrix}0&\elltilde_i&0\end{smallmatrix}\right)$ and $r_i=\left(\begin{smallmatrix}0\\\rtilde_i\\0\end{smallmatrix}\right)$ then $\Pi_d=\sum_ir_i\cdot\ell_i$. Finally, $F^i\cdot r_t=\rho^i r_{t+i}$ where $t+i$ is taken modulo~$d$. Henceforth:
\begin{gather*}
F^i\cdot  \Pi_d =\rho^i\sum_t r_{t+i}\cdot\ell_{t}
\end{gather*}
and thus, according to Th.~\ref{thr:3}:
\begin{align*}
  \Theta&=\frac1d\Bigl(\,\sum_{i=0}^{d-1}\rho^{-i}F^i\Bigr)\cdot\Pi_d\\
        &=\frac1d\Bigl(\,\sum_{i=0}^{d-1}\sum_{t=0}^{d-1}r_{t+i}\cdot\ell_t\Bigr)\\
  &=\frac1d\bigl(\,\sum_ir_i\Bigr)\cdot\bigl(\sum_t\ell_t\bigr)
\end{align*}
But $r=\frac1d\sum_i r_i$ and $\ell=\sum_i\ell_i$ form a pair of right and left Perron eigenvectors of $F$ satisfying $\ell\cdot r=\frac1d\sum_{i,j}\delta_i^j=1$, whence the announced result.

\subsection{Proof of Theorem~\ref{thr:6}}

Assume that $W$ is an umbrella digraph, and let $\rho=\rho(W)$. It follows from Th.~\ref{thr:3} that the adjacency matrix $F$ has a positive $\rho$-eigenvector, say~$r$. Consider the positive cocyle $\Gamma(x,y)=\evm ry/\evm rx$ and the following cocycle kernel:
\begin{gather}
  \label{eq:31}
  q(x,y)=\rho^{-1}\w(x,y)\Gamma(x,y)
\end{gather}
This is indeed a transition kernel since, for every~$x$:
\begin{align*}
  \sum_yq(x,y)=\frac{\rho^{-1}}{\evm rx}\sum_y \evm F{x,y}\evm ry=\frac{\rho^{-1}}{\evm rx}\evm{F\cdot r}x=1
\end{align*}
Hence $q$ defined in~(\ref{eq:31}) is a cocycle kernel which defines a complete cocycle measure.

Furthermore, if $W$ is accessible from a given vertex~$x_0$, one shows that all the complete cocycle measures are of the form~(\ref{eq:31}), where the positive cocycles $\Gamma$ are in bijection with positive $\rho$-eigenvectors~$r$, up to a multiplicative constant, through:
\begin{gather*}
  \evm rx=\Gamma(x,x_0)
\end{gather*}
Since positive $\rho$-eigenvectors are exactly of the form $\sum_i\alpha_ir_i$ for all families $(\alpha_i)_i$ of positive reals according to Th.~\ref{thr:3}, the parameterizing of complete cocycle measures by an open simplex of dimension $p-1$ follows.

\medskip
Conversely, we assume that $W$ has a complete cocycle measure $\mu=(\mu_x)_{x\in V}$. Let us prove that $W$ is necessarily an umbrella digraph. Let $q(x,y)=\lambda\w(x,y)\Gamma(x,y)$ be the cocycle kernel associated with~$\mu$. Since the cocycle measure $\mu$ is assumed to be complete, we fix a positive bound $m>0$ of $\Gamma$ on~$\Tos$.

For every vertex~$y$, let $\rho(y)$ denote the spectral radius of the access class of~$y$. Let $x$ be a vertex. The spectral radius of $V(x)$ is $\gamma(x)=\max_{y\in V(x)}\rho(y)$. Let $y\in V(x)$ be a vertex such that $\rho(y)=\gamma(x)$. For brevity, we put $\rho=\rho(y)$. Then $\rho^{-1}$ is the convergence radius of the power series $\sum_k Z_y(k)z^k$. Pick a path $u\in O_{x,y}$. Then, for every integer $k\geq0$:
  \begin{gather}
    \sum_{v\in O_y(k)}\mu_x\bigl(\up(uv)\bigr)\leq 1
  \end{gather}
  which yields $Z_y(k)K\leq\lambda^{-k}$ with $K=\lambda^{\length u}\w(u)m$, and therefore $\rho^{-1}\geq\lambda$.

  Let $F$ be the adjacency matrix of~$V(y)$. Then we claim that $\lambda^{-1}$ is an eigenvalue of~$F$. Let $\alpha=(\alpha_p)_{p\in V(y)}$ be defined by $\alpha(p)=\Gamma(y,p)$ for $p\in V(y)$. For every $p\in V(y)$, one has:
  \begin{gather}
\label{eq:33}    
    \sum_{q\in V(y)}\mu_p(\up q)=1,\quad\text{yielding:\quad}
    \sum_{q\in V\tq (p,q)\in E}\lambda\w(p,q)\Gamma(p,q)=1
  \end{gather}
  Every $q\in V$ such that $(p,q)\in E$ satisfies $y\To p\To q$ and thus $\Gamma(y,q)=\Gamma(y,p)\Gamma(p,q)$. All the terms in this equation are positive, hence $\Gamma(p,q)=\alpha(q)/\alpha(p)$, and thus~\eqref{eq:33} yields:
  \begin{gather}
    \lambda^{-1}\alpha(p)=\sum_{q\in V\tq (p,q)\in E}\w(p,q)\alpha(q)
  \end{gather}
  proving that $\lambda^{-1}$ is an eigenvalue of~$F$. Since $\rho\leq\lambda^{-1}$ and since $\rho$ is the largest eigenvalue of~$F$, it implies that $\rho=\lambda^{-1}$, hence $\gamma(x)=\lambda^{-1}$ for every $x\in V$. But $\rho(W)=\max_{x\in V}\gamma(x)$, so $\rho(W)=\lambda^{-1}$.

  Let $x$ be a vertex in some final access equivalence class of~$W$. Then $\mu$ induces a complete cocycle measure on the access class~$D(x)$ of~$x$. Therefore, by the previous point, $\rho(x)=\gamma(x)$ and thus $\rho(x)=\lambda^{-1}$. This proves that every final access equivalence is a basic class.

  Now let $x\in V$ be a vertex such that $\rho(x)=\lambda^{-1}$. We claim that $x$ belongs to a final access equivalence class.

  Seeking a contradiction, assume that it is not.  Let $F$ be the adjacency matrix of~$V(x)$. After maybe a simultaneous permutation of the lines and the columns, we may assume that $F$ has the following form by blocks:
\begin{gather}
  F=
  \begin{pmatrix}
    A&B\\
    0&K
  \end{pmatrix}
\end{gather}
where $A$ is the adjacency matrix of the access class~$D(x)$ of~$x$. The matrix $A$ is irreducible. Let $\alpha$ be defined, as above, by $\alpha(p)=\Gamma(x,p)$ for $p\in V(x)$, and let $\beta$ be the restriction of $\alpha$ to~$D(x)$. Then, since $B$ is non negative and non zero on the one hand, and since $F\alpha=\lambda^{-1}\alpha$ and $\alpha>0$ on the other hand, one gets $A\beta\leq \lambda^{-1}\beta$ with $A\beta\neq \lambda^{-1}\beta$. But, since $\beta>0$, this contradicts the subinvariance theorem for irreducible matrices~\cite[Th.~1.6]{seneta81}, which proves our claim.

\subsection{Proof of Corollary~\ref{cor:1}}
\label{sec:proof-corollary}

It derives from the correspondence between positive $\rho$-eigenvectors and complete cocycle measures, which are proved in Th.~\ref{thr:6} to exist for $W$ if and only if $W$ is an umbrella digraph.

\section{Proofs for Section~\ref{sec:comp-asympt-matr}}
\label{sec:proofs-section}

\subsection{Proof of Lemma~\ref{lem:2}}
\label{sec:proof-lemma}

  If $F_S$ and $F_T$ denote the adjacency matrices of $S$ and~$T$, one has, for each integer $n\geq0$:
  \begin{align*}
    F^n&=
    \begin{pmatrix}
      F_S^n&Y_n\\
      0&F_T^n
    \end{pmatrix}
         &Y_n&=\sum_{k=0}^{n-1}F_S^{k}\cdot X\cdot F_T^{n-1-k}
  \end{align*}
yielding for $(x,y)\in S\times T$\,:
\begin{align*}
\evM { Y(z)}{x,y}&=z\sum_{(u,v)\in S\times T}\evM X{u,v}
\sum_{n\geq1}z^{n-1}\Bigl(\,\sum_{k=0}^{n-1}\evM{F_S^k}{x,u}\evM{F_T^{n-1-k}}{v,y}\Bigr)\\
                &=z\sum_{(u,v)\in S\times T}
                  \evM X{u,v}\evM{\H_S(z)}{x,u}\evM{\H_T(z)}{v,y}\\
  &=z\evM{\H_S(z)\cdot X\cdot\H_T(z)}{x,y}
\end{align*}
which proves the formula~(\ref{eq:4}).

\subsection{Proof Theorem~\ref{thr:1}}
\label{sec:proof-theorem-1}

The proof can be done by induction on the number of access classes of~$W$, yielding the same discussion as in the core of the text, Section~\ref{sec:recursive-computing}, \emph{Recursive computing}. Some further discussion is needed to recursively determine the pairs $(x,y)$ such that $\evm\Theta{x,y}>0$, without additional difficulty.

\section{Proofs for Section~\ref{sec:conv-distr}}

\subsection{Proof of Theorem~\ref{thr:9}}
\label{sec:proof-theorem-2}

The proof uses the same computations as the ones performed in the core of the text, Section~\ref{sec:itsh-unif-distr}, \emph{Uniform distributions}. It is completed by the form of the residual matrix of $W$ given in Corollary~\ref{cor:2}.

\subsection{Proof of Corollary~\ref{cor:3}}
\label{sec:proof-corolalry}

The result is clear if $\rho=0$, hence we assume that $\rho>0$ for the rest of the proof.
Let $\varphi(X)$ be the minimal polynomial of~$F$~:
\begin{gather*}
  \varphi(X)=(X-\lambda_1)^{m_1}(X-\lambda_2)^{m_2}\dots(X-\lambda_q)^{m_q}
\end{gather*}
where $\lambda_1=\rho$ and $\lambda_2,\dots,\lambda_q$ are the distinct eigenvalues of $F$ and $m_1,\ldots,m_q$ are their multiplicities in the minimal polyonomial. Hence $m_1,\ldots,m_q$ are also the dimensions of the associated generalized eigenspaces.

It follows from the theory of functions of matrices (see F.R.~Gantmacher, \emph{The Theory of Matrices}, Vol.~1, Chap.~V) that there are $m$ square matrices $T_{i,j}$ for $i=1,\ldots,q$ and $j=1,\ldots,m_i$, where $m=m_1+\dots+m_q$, which are called the \emph{components of the matrix~$F$}, such that for every complex polynomial~$f$, one has:
\begin{gather}
  \label{eq:34}
  f(F)=\sum_{i=1}^q\Bigl(\,\sum_{j=0}^{m_i-1} f^{(j)}(\lambda_i)T_{i,j}\Bigr)
\end{gather}

Applying~\eqref{eq:34} in particular with the polynomial $f_k(z)=z^k$ for every nonnegative integer~$k$, and then summing over $k$ yields:
\begin{align*}
  \H(z)&=
\sum_{i,j}j!\frac{(\lambda_i z)^j}{(1-\lambda_iz)^{j+1}}T_{i,j}\\
  (1-\rho s)^h\H(s)&=\sum_{j=0}^{m_1-1}j!(\rho s)^j(1-\rho s)^{h-j-1}T_{1,j}
                   +\sum_{i\neq 1,j}j!\frac{(\lambda_i z)^j(1-\rho s)^h}{(1-\lambda_iz)^{j+1}}T_{i,j}
\end{align*}

The components $T_{i,j}$ are linealry independent (\emph{cf.} Gantmacher, \emph{ibidem}), and in particular none of them is zero. Therefore if $m_1>h$, the above expression shows that $(1-\rho s)^h\,\H(s)$ cannot have a finite limit at~$\rho^{-1}$, contrary to the statement of Theorem~\ref{thr:1}. If $m_1<h$ then the above expressions shows that $(1-\rho s)^h\,\H(s)$ has limit zero at~$\rho^{-1}$, contradicting again the statement of Th.~\ref{thr:1}. Hence $m_1=h$, which was to be shown.

\subsection{Proof of Theorem~\ref{thr:10}}

This is done by studying the convergence of each subsequence $(\mu_{x,kd+i})_{k\geq0}$ for $i\in\{0,\ldots,d-1\}$. The sequence $(\mu_{x,k})_{k\geq0}$ converges if and only if all the subsequences converge toward the same limit, which is seen to be equivalent to the given condition.

\end{document}